\documentclass[preprint,12pt]{elsarticle}

\usepackage{amsmath,amssymb,mathtools}
\usepackage{algorithm}
\usepackage{algpseudocode}
\usepackage{graphicx}
\usepackage{booktabs}
\usepackage{hyperref}

\hypersetup{
  colorlinks=true,
  linkcolor=blue,
  citecolor=blue,
  urlcolor=blue
}

\newcommand{\R}{\mathbb{R}}

\newcommand{\T}{\mathbb{T}}

\newcommand{\Fcal}{\mathcal{F}}
\newcommand{\grad}{\nabla}

\journal{Applied Mathematics Letters}

\begin{document}

\begin{frontmatter}

\title{Fourier-Diagonalized Natural Gradients and Sobolev Mirror Descent}

\author[1]{Jeongbin Jo\corref{cor1}}
\ead{jeongbin033@yonsei.ac.kr}
\cortext[cor1]{Corresponding author}
\address[1]{Yonsei University, Seoul, Republic of Korea}

\begin{abstract}
We study natural-gradient updates whose metric operators are diagonalized by the Fourier transform and relate them to Sobolev mirror descent. Translation-invariant Fisher geometries and Sobolev mirror geometries share a common inverse-map structure in the spectral domain. The Fisher metric is represented by a positive Fourier symbol, while Sobolev mirror geometry corresponds to the specific Bessel-potential symbol associated with the Sobolev norm. When these symbols coincide, the natural-gradient and mirror-descent updates are identical; otherwise, Sobolev mirror descent provides a canonical spectral preconditioner for the Fisher inverse geometry. This gives a mathematical lens through which spectral filtering and truncation techniques in PDE and operator learning can be viewed as natural actions of inverse metric geometry. We introduce Spectral Natural Gradient, an FFT-based implementation of these geometric updates.
\end{abstract}


\end{frontmatter}

\section{Introduction}
\label{sec:intro}

Optimizing parameterized probability distributions and continuous spatial fields is a central problem across machine learning, physics-informed neural networks, and statistical mechanics. Ordinary first-order gradient descent assumes the parameter space is flat and Euclidean, often leading to slow convergence or unstable behavior in ill-conditioned landscapes, especially when oscillatory spatial modes are present~\citep{nocedal2006numerical,wang2021understanding}.

To account for the intrinsic curvature of these spaces, Amari introduced Natural Gradient Descent (NGD)~\citep{amari1998natural, amari2016information, martens2020new}. NGD operates within the framework of information geometry, utilizing the Fisher Information Matrix (FIM), denoted as $F(\theta)$, as the Riemannian metric tensor. The exact NGD update rule requires evaluating $\theta_{t+1} = \theta_t - \alpha\, F^{-1}(\theta_t)\,\grad L(\theta_t)$. By following the steepest descent direction in distribution space rather than Euclidean parameter space, NGD provides a geometrically meaningful scaling of the gradient. However, for a model with $N$ parameters, explicitly forming a dense FIM requires $\mathcal{O}(N^2)$ storage, and direct inversion or factorization can require $\mathcal{O}(N^3)$ operations. This cubic scaling creates a formidable computational bottleneck for extensive spatial models.

In this work, we establish a spectral dictionary between translation-invariant natural-gradient geometry and Sobolev mirror descent. Both frameworks share a common Fourier-diagonal inverse-map structure. While the translation-invariant Fisher metric takes an arbitrary positive Fourier symbol $\lambda_F(\omega)$, Sobolev mirror descent is the special case in which the spectral symbol is $(1+|\omega|^2)^s$. When these symbols coincide, the two updates are identical; otherwise, Sobolev mirror descent provides a canonical spectral preconditioner for the Fisher inverse geometry. This dictionary gives an information-geometric explanation for spectral filtering mechanisms used in spatial optimization.

This perspective is complementary to scalable natural-gradient and preconditioned optimization methods that impose tractable matrix structure, including truncated or approximate natural-gradient solvers~\citep{pascanu2014revisiting}, Kronecker-factored curvature approximations~\citep{martens2015kfac}, and tensor-structured preconditioners~\citep{gupta2018shampoo}. Rather than approximating a generic dense Fisher matrix, SNG isolates the case in which the relevant geometry is already diagonal in a known harmonic basis. It is also related to Sobolev-gradient methods for differential equations~\citep{neuberger1997sobolev}, but emphasizes the natural-gradient interpretation and the explicit FFT implementation of the inverse geometry.

\section{Fourier-Diagonalized Geometry and Sobolev Mirror Descent}
\label{sec:duality}

In this section, we formulate the geometry at the level of continuous function spaces. We first state the spectral linear-algebra fact behind the inverse map, and then use mirror descent to interpret one important instance: the Sobolev geometry generated by $(I-\Delta)^s$. The finite-dimensional FFT algorithm introduced later should be understood as a periodic discretization of this continuous Fourier-multiplier structure.

\paragraph{Domain scope and spectral direction.}
We consider fields on either the whole space $\R^d$ or the $d$-dimensional torus $\T^d$. In these settings, translation-invariant operators are Fourier multipliers: on $\R^d$ this is expressed by the continuous Fourier transform, while on $\T^d$ it is expressed by Fourier series. The frequency variable $\omega$ denotes the corresponding wavenumber, scaled by the chosen length unit when the physical domain has dimensions. General bounded domains with non-periodic boundary conditions are outside the exact Fourier-multiplier setting considered here.

Let $F$ be a positive self-adjoint translation-invariant operator on $\R^d$ or $\T^d$, possessing a positive Fourier symbol $\lambda(\omega) > 0$. The natural-gradient direction $d = F^{-1}g$ is given modewise by
\begin{equation}
  \widehat{d}(\omega) = \frac{\widehat{g}(\omega)}{\lambda(\omega)} .
  \label{eq:spectral_ng_direction}
\end{equation}
Because translation invariance implies that $F$ is diagonal in the Fourier representation, we have $\widehat{Fd}(\omega) = \lambda(\omega)\widehat{d}(\omega)$. Solving $Fd = g$ therefore reduces to scalar division for each frequency. In a periodic finite-dimensional discretization, the same statement becomes diagonalization by the discrete Fourier transform and can be implemented by FFT.

\subsection{Mathematical Definition of Mirror Descent (Fenchel Duality)}
\label{sec:mirror_def}

We begin by defining the mechanics of Mirror Descent~\citep{nemirovski1983problem, beck2003mirror}. Mirror descent generalizes gradient descent beyond Euclidean settings and is closely related to information-geometric optimization~\citep{raskutti2015information}. Let the primal space be $\mathcal{X}$, and its dual space be $\mathcal{X}^*$. Given a strictly convex distance-generating function $\Phi(x)$, the \emph{Mirror Map} from the primal space to the dual space is defined as the gradient of $\Phi$:
\begin{equation}
  y = \nabla \Phi(x) \quad (x \in \mathcal{X}, \; y \in \mathcal{X}^*).
\end{equation}

The conjugate of this function, known as the Fenchel Conjugate $\Phi^*(y)$, is defined via the Legendre-Fenchel transform:
\begin{equation}
  \Phi^*(y) = \sup_{x \in \mathcal{X}} \left( \langle x, y \rangle - \Phi(x) \right).
\end{equation}

By standard convex duality, the map returning from the dual space back to the primal space---the \emph{Inverse Mirror Map}---is the gradient of the Fenchel conjugate:
\begin{equation}
  x = \nabla \Phi^*(y) = (\nabla \Phi)^{-1}(y).
\end{equation}

Therefore, the $t$-th iteration of Mirror Descent is a composition of two distinct steps:
\begin{enumerate}
    \item \textbf{Gradient descent in the dual space:} $y_{t+1} = \nabla \Phi(x_t) - \alpha \nabla J(x_t)$
    \item \textbf{Return to the primal space:} $x_{t+1} = \nabla \Phi^*(y_{t+1})$
\end{enumerate}

\subsection{Extension to Infinite-Dimensional Function Spaces (Sobolev Space)}
\label{sec:sobolev}

We now expand the variable $x$ from a finite array of numbers to a continuous function $u(v)$ parameterized by spatial coordinates $v$. The objective function then becomes a functional $J(u)$ on a function space.

The key modeling choice is the distance-generating function $\Phi(u)$. To control the curvature and smoothness of the space, we replace the standard $L_2$ Euclidean norm with the Sobolev $H^s$ norm~\citep{brezis2010functional}. In the Fourier representation, this norm is written uniformly as
\begin{equation}
  \Phi(u) = \frac{1}{2} \| u \|_{H^s}^2
  = \frac{1}{2}\int_{\widehat{\Omega}}
  (1+|\omega|^2)^s |\hat u(\omega)|^2\,d\mu(\omega),
\end{equation}
where $\widehat{\Omega}=\R^d$ with Lebesgue measure for the whole-space case, and $\widehat{\Omega}=\mathbb{Z}^d$ with counting measure for the torus case.

Calculating the Fréchet derivative of $\Phi(u)$ yields a linear elliptic operator:
\begin{equation}
  \nabla \Phi(u) = (I - \Delta)^s u,
\end{equation}
where $I$ is the identity operator and $\Delta$ is the Laplacian. This operator is a Bessel-potential form of the fractional Laplacian, a standard object in fractional Sobolev analysis~\citep{dinezza2012hitchhiker}. The mirror map $\nabla \Phi$ weights higher spatial frequencies more strongly, while the inverse mirror map suppresses them.

\subsection{Exact Implementation via Fourier Transform}
\label{sec:fourier_impl}

We substitute this Sobolev mirror map into the generic Mirror Descent update equations.

\paragraph{Step 1: Update in Dual Space.}
\begin{equation}
  \label{eq:dual_update}
  y_{t+1} = (I - \Delta)^s u_t - \alpha \nabla J(u_t).
\end{equation}
Here, $y_{t+1}$ is an element of the dual space containing the high-order derivative information of the function $u$.

\paragraph{Step 2: Return to Primal Space via Inverse Mirror Map.}
To obtain the new function $u_{t+1}$ in the primal space, we must apply the inverse differential operator:
\begin{equation}
  u_{t+1} = (\nabla \Phi)^{-1} y_{t+1} = (I - \Delta)^{-s} y_{t+1}.
\end{equation}

The Fourier representation now becomes decisive. Computing the inverse of the differential operator $(I - \Delta)^{-s}$ directly in the spatial domain is expensive if treated as a dense linear solve. However, spectral methods~\citep{trefethen2000spectral} exploit the fact that the spatial Laplacian $\Delta$ is represented by scalar multiplication by $-|\omega|^2$ in frequency variables:
\begin{equation}
  \mathcal{F}[\Delta u](\omega) = -|\omega|^2 \mathcal{F}[u](\omega).
\end{equation}

Consequently, the inverse of this differential operator becomes scalar division in the frequency domain:
\begin{enumerate}
    \item Take the Fourier representation of the dual variable $y_{t+1}$:
    \begin{equation}
        \hat{y}_{t+1}(\omega) = \mathcal{F}[y_{t+1}].
    \end{equation}
    \item Perform the inverse mirror map operation as a simple division in frequency space:
    \begin{equation}
        \hat{u}_{t+1}(\omega) = \frac{1}{(1 + |\omega|^2)^s} \hat{y}_{t+1}(\omega).
    \end{equation}
    \item Transform back to the primal function $u_{t+1}(v)$:
    \begin{equation}
        u_{t+1}(v) = \mathcal{F}^{-1} \left[ \frac{1}{(1 + |\omega|^2)^s} \mathcal{F}[y_{t+1}] \right].
    \end{equation}
\end{enumerate}

Equivalently, composing the dual update with the inverse mirror map gives the primal update
\begin{equation}
  u_{t+1}
  = u_t - \alpha (I-\Delta)^{-s}\nabla J(u_t),
\end{equation}
which is precisely a spectral preconditioned gradient step with symbol
$\lambda(\omega)=(1+|\omega|^2)^s$.

\subsection{Mathematical Conclusion and Physical Meaning}
\label{sec:conclusion_physical}

Through this derivation, we arrive at the following conclusion for Sobolev mirror descent:

\begin{quote}
\textit{"The inverse mirror map generated by the Sobolev $H^s$ norm is a low-pass spectral filter in the frequency domain."}
\end{quote}

Consider the analytical form of the scalar multiplier: $\frac{1}{(1 + |\omega|^2)^s}$.
\begin{itemize}
    \item As $\omega$ increases, the denominator grows and high-frequency modes are damped.
    \item As $\omega$ decreases, the multiplier approaches $1$ and preserves low-frequency structure.
\end{itemize}

Thus, Sobolev mirror descent embeds a mathematically explicit Fourier spectral filter directly into the optimization step. This establishes our spectral dictionary: it shows why spectral low-pass filtering is natural in PDE and operator-learning optimization. From this perspective, the spectral truncation employed in architectures such as Fourier Neural Operators (FNO)~\citep{li2021fourier} can be viewed as a band-limited analogue of Sobolev inverse filtering, suggesting a geometric interpretation of spectral truncation in operator-learning models.

\section{Spectral Natural Gradient (SNG) Algorithm}
\label{sec:sng}

We now pass from the continuous Fourier-multiplier update to its finite-dimensional periodic implementation. The \textbf{Spectral Natural Gradient (SNG)} computes natural-gradient or mirror-descent updates whenever the relevant geometry has a known Fourier symbol $\lambda(\omega)$. It is exact for Fourier-diagonalizable Fisher operators and for Sobolev mirror-descent metrics, while serving as a structured spectral preconditioner when $\lambda(\omega)$ is used as an approximation.

\begin{algorithm}[ht]
  \caption{Spectral Natural Gradient (SNG)}
  \label{alg:sng}
  \begin{algorithmic}[1]
    \Require parameters $u_t$, step size $\alpha$, spectral geometry $\lambda(\omega)$, regularizer $\epsilon>0$
    \State \textbf{Primal gradient computation:}\;
      $g_t \leftarrow \grad J(u_t)$
    \State \textbf{Fourier transform:}\;
      $\widehat{g}_t(\omega) \leftarrow \Fcal[g_t](\omega)$
    \State \textbf{Spectral inverse geometry via scalar division:}\;
      $\widehat{d}_t(\omega) \leftarrow \dfrac{\widehat{g}_t(\omega)}{\lambda(\omega)+\epsilon}$
    \State \textbf{Inverse transform and parameter update:}\;
      $u_{t+1} \leftarrow u_t - \alpha\,\Fcal^{-1}[\widehat{d}_t]$
  \end{algorithmic}
\end{algorithm}

For a Sobolev geometry, $\lambda(\omega)=(1+|\omega|^2)^s$. For a shift-invariant Fisher or Hessian operator, $\lambda(\omega)$ is the Fourier symbol of that operator. After periodic discretization on $N$ grid points, the total computational cost per step is bounded by the Fast Fourier Transforms and pointwise divisions, giving $\mathcal{O}(N \log N)$ time.

\section{Numerical Experiments}
\label{sec:experiments}

The preceding sections establish the Sobolev filtering interpretation analytically in function space. We next give a numerical illustration of the exact spectral natural-gradient case after periodic discretization. We construct a shift-invariant 1D quadratic convolution objective, $L(u)=\tfrac12\|x*u-y\|^2$. This objective may be viewed as the negative log-likelihood of a Gaussian linear convolution model with fixed noise variance; in this setting the Fisher operator and the quadratic Hessian coincide up to the noise scale. Because the discrete convolution is periodic, the resulting Fisher/Hessian matrix is circulant. We evaluate both algorithms: exact NGD using dense linear algebra, and SNG using FFT and scalar division by the exact circulant spectrum.

\begin{figure}[ht]
  \centering
  \includegraphics[width=\textwidth]{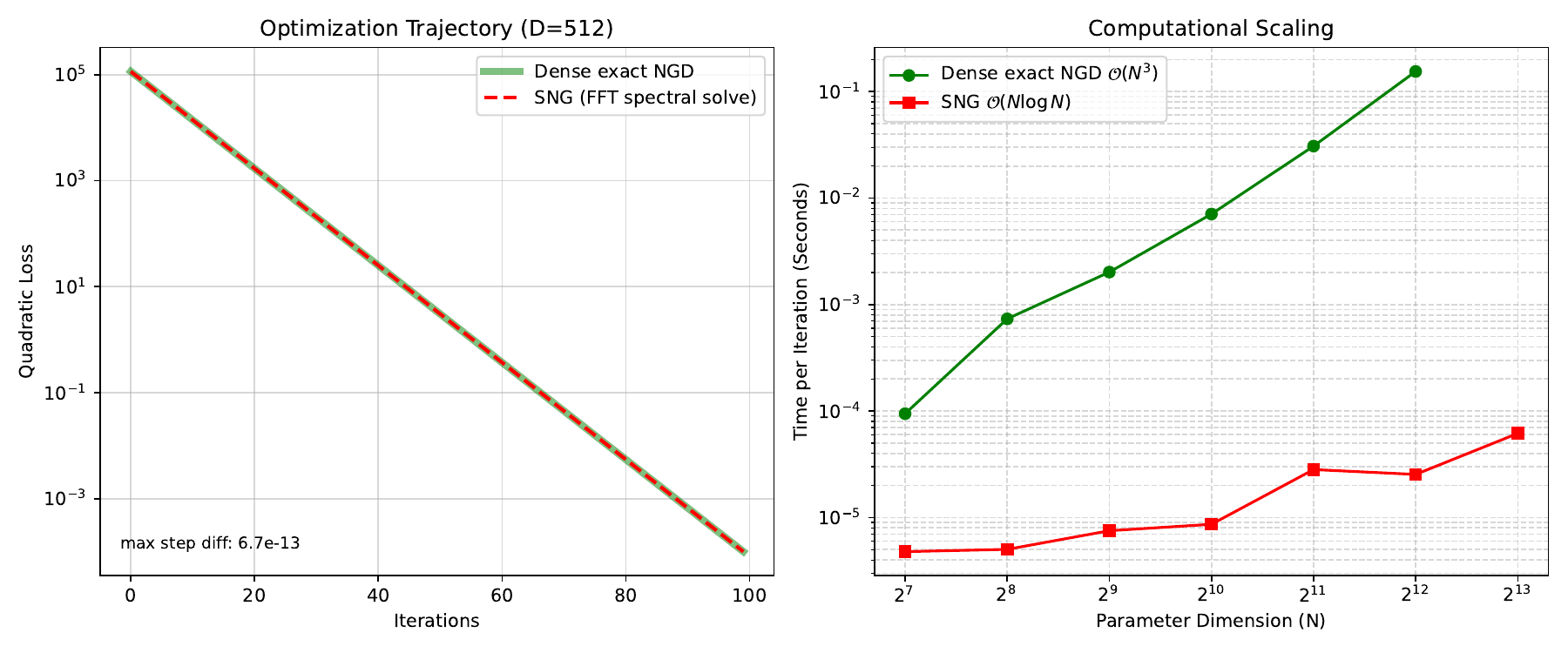}
  \caption{Comparison of Exact NGD and SNG on a shift-invariant quadratic convolution problem. \textbf{Left:} Optimization trajectory for $N=512$; the overlap verifies equivalence between the dense circulant solve and Fourier-domain scalar division. \textbf{Right:} Log-log plot of time per iteration vs. parameter dimension $N$, illustrating the $\mathcal{O}(N^3)$ vs. $\mathcal{O}(N \log N)$ scaling advantage.}
  \label{fig:time_comparison}
\end{figure}

Figure~\ref{fig:time_comparison} (left) plots the optimization trajectories for a parameter dimension of $N=512$. The loss curves for SNG and exact NGD overlap, confirming that Fourier-domain scalar division implements the same linear solve as dense Fisher/Hessian inversion under shift-invariant constraints. In double precision, the maximum step-wise discrepancy between the dense and spectral directions was $6.7\times 10^{-13}$.

Figure~\ref{fig:time_comparison} (right) compares per-step runtime as the parameter dimension $N$ ranges from $128$ to $8192$. The timings use double precision and averages over 10 repeated solves; their purpose is to illustrate scaling rather than hardware-specific performance. Dense exact NGD follows the expected cubic scaling of direct linear solves, while SNG follows the FFT-based $\mathcal{O}(N \log N)$ scaling and avoids explicit construction or inversion of the dense circulant matrix.

\section{Conclusion}
\label{sec:conclusion}

Natural gradient descent is a geometric optimizer traditionally hampered by dense matrix storage and cubic inversion or factorization costs. We showed that for Fourier-diagonalizable geometries, including translation-invariant Fisher/Hessian operators and Sobolev mirror-descent metrics, this dense inverse can be replaced by FFT-based scalar division in the spectral domain. In the Sobolev case, the inverse mirror map reduces precisely to a low-pass multiplier. Our proposed Spectral Natural Gradient (SNG) therefore achieves exact natural-gradient trajectories for shift-invariant Fisher operators and exact Sobolev mirror-descent trajectories for Laplacian metrics, both in $\mathcal{O}(N \log N)$ time. The exactness claim is restricted to metrics with a known Fourier symbol; outside this class, SNG should be interpreted as a structured spectral preconditioner. This provides a precise route for applying geometric optimization ideas to PDE and spatial learning domains where spectral structure is available.

\bibliographystyle{elsarticle-num}
\bibliography{references}

\end{document}